\documentstyle[12pt]{article}

\begin{document}

 \def\bide{\mbox{$\Delta$}}
 \def\bidet{\mbox{$\widetilde{\Delta}$}}
 \def\ggt{\mbox{$\widetilde{g}$}}
 \def\SSt{\mbox{$\widetilde{S}$}}
 \def\rrt{\mbox{$\widetilde{\rho}$}}
 \def\epw{\mbox{$\epsilon(w)$}}
 \def\epv{\mbox{$\epsilon(v)$}}
 \def\CCt{\mbox{$\widetilde{C}$}}
 \def\ovC{\mbox{$\overline{C}$}}
 \def\ovCCt{\mbox{$\overline{\CCt}$}}
 \def\carla{\mbox{$ch \rule{1mm}{0mm} L^{\lambda}$}}
 \def\carmu{\mbox{$ch \rule{1mm}{0mm} L^{\mu}$}}
 \def\carks{\mbox{$ch \rule{1mm}{0mm} L^{\xi}$}}
 \def\lla{\mbox{$L^{\lambda}$}}
 \def\psil{\mbox{$\Psi^{\lambda}$}}
 \def\psim{\mbox{$\Psi^{\mu}$}}
 \def\psitm{\mbox{$\widetilde{\Psi^{\mu}}$}}
 \def\nmti{\mbox{$\widetilde{n_{\mu}}$}}
 \def\nnti{\mbox{$\widetilde{n_{\nu}}$}}
 \def\nxti{\mbox{$\widetilde{n_{\xi}}$}}
 \def\proal{\mbox{$\prod_{\alpha \in \Delta} (1- e^{-\alpha})$}}
 \def\promi{\mbox{$\prod_{\Delta \setminus \bidet}
 (1- e^{-\alpha})$}}
 \def\Kgtg{\mbox{$K_{\textstyle \ggt \subset g}$}}
\thispagestyle{empty}
\begin{flushright}
{\large SPBU-IP-95-8}
\end{flushright}
\vspace*{2cm}

\begin{center}
{\bf \large
Recursion relations and branching rules \\
for simple Lie algebras}\\[1cm]

{\bf
 V.D.Lyakhovsky
 \footnote{ E-mail address: LYAKHOVSKY @ PHIM.NIIF.SPB.SU } \\
 S.Yu.Melnikov \\
 Theoretical Department \\
 Institute of Physics \\
 St.Petersburg State University\\
 St.Petersburg \\
 198904, Russia
 } \\[1cm]
  {\bf Abstract.}
\end{center}

The branching rules between simple Lie algebras and its regular
(maximal) simple subalgebras are studied. Two types of recursion
relations for anomalous relative multiplicities are obtained.
One of them is proved to be the factorized version of the other.
The factorization property is based on the existence of the set
of weights $\Gamma$ specific for each injection. The structure of
$\Gamma$ is easily deduced from the correspondence between the
root systems of algebra and subalgebra. The recursion relations thus
obtained give rise to simple and effective algorithm for branching
rules. The details are exposed by performing the explicit
decomposition procedure for $ A_{3} \oplus u(1) \rightarrow B_{4} $
injection.

\newpage
\setcounter{page}{1}

\section{Introduction}
\subsection{}
In elementary particle physics and especially in model building it
is quite important to have effective branching rules for Lie algebras
representations. There are several simple methods of decomposition
appropriate for different types of injections, for example,
the Gelfand-Zetlin method for $A_{n-1} \rightarrow A_n ,
B_{n-1} \rightarrow B_n $ and
$D_{n-1} \rightarrow D_n $ \cite{G-Z,L} .
In the general case the most advanced investigation was
performed by R.V.Moody,
J.Patera,R.T.Sharp and F.Gingras in a series of works
\cite{Pat,King,Sha}.
Their approach is
based on the properties of Weyl orbits in weight diagrams and the
generating function technique \cite{GF}.

In this paper we want to demonstrate that recursion relations for
multiplicities of subrepresentations can also be successfully used
in the decomposition procedure. We restrict the exposition to regular
maximal injections of reductive subalgebras (composed of semisimple
ones and diagonalizable Abelian algebras), though it is possible to
treat analagously special injections and also nonmaximal ones.
For a simple algebra $g$ and its regular maximal reductive
subalgebra $\ggt$ the problem is to evaluate the coefficients
$n_{\mu}$ in the decomposition of an irreducible representation
$\lla(g)$ ($\lambda$ is its highest weight)
\begin{equation}
  \lla(g)_{\downarrow \ggt} =
  \oplus_{\mu} n_{\mu} \widetilde{L}^{\mu}(\ggt).   \label{e1}
\end{equation}
For any weight $\nu$ of $\lla$ the total multiplicity
$m_{\nu}$ can be presented as the sum
\begin{equation}
m_{\nu} = m_{\nu}' + n_{\nu}               \label{e2}
\end{equation}
where $m_{\nu}'$ is the multiplicity induced by the
subrepresentation $\oplus_{\mu > \nu} n_{\mu}
\widetilde{L}^{\mu}(\ggt)$ contained in (\ref{e1}).
The second term $n_{\nu}$ in (\ref{e2}) is called {\em the relative
multiplicity} of the weight $\nu$. When $\nu$ is from the dominant
Weyl chamber its relative multiplicity coinsides with the
corresponding coefficient in the decomposition (\ref{e1}).
The recursion relations for relative multiplicities for the injection
$ A_3 \oplus u(1) \rightarrow D_4 $ were studied in \cite{L-F}.
It was shown that the considerable set of multiplicities for
intermediate weights mutually cancel and the final recursion formula
is suitable for calculations. As will be demonstrated below the
general recursion relation for regular injections can be naturally
formulated in terms of {\em anomalous relative multiplicities}
$\nnti$. To define them consider the highest weights
$$
M = \{ \mu \mid n_{\mu} \neq   0 \}
$$
for the decomposition (\ref{e1}). For the subalgebra $\ggt$
let $V$ be the Weyl group and $\rrt$ --
the half-sum of positive roots. The anomalous relative
multiplicity \nnti $(g,\ggt,\lambda,\nu)$
is the function
\begin{equation}
\widetilde{n} (g,\ggt,\lambda,\nu) \equiv
\nnti = \left\{ \begin{array}{ll}
                det(v)n_{\mu} & \mbox{for $ \{ \mu \in M \mid
                v(\mu + \rrt) - \rrt = \nu  \} $}\\
                0             & \mbox{elsewhere}
                \end{array}
        \right.
\label{e3}
\end{equation}
defined on the weight space of $g$.
As will be demonstrated the general
recursion relation for regular injections can be naturally
formulated in terms of
anomalous relative multiplicities \nnti.

The paper is organized as follows. The general formalism is presented
 in section 2. The obtained recursion relations are based on the
 properties of the elementary "fan" $\Gamma$ -- the special set of
 weights defined by the injection $ \ggt  \rightarrow g $. The
 structure of $\Gamma$'s and their basic properties are studied
in details. To demonstrate explicitly the role of
$\Gamma$ in recursion procedure
we use the very simple example -- the injection
$ A_1 \oplus u(1) \rightarrow B_2 $. In the appendix the
applications of the recursion formulas are shown in full details for
the injection $ A_3 \oplus u(1) \rightarrow B_4  $.

\subsection{Notations}

$\rule{6mm}{0mm}g$ -- the simple Lie algebra;

\ggt \, -- the reductive regular subalgebra of $g$;

\bide, \bidet \, -- the corresponding sets of positive roots, note
that \bidet \, is the system of positive roots of the semisimple
subalgebra in \ggt ;

$S$, \SSt \, -- the sets of basic roots;

$\rho$, \rrt \, -- the half-sums of positive roots for $g$ and \ggt
\, respectively;

$W,V$ \, -- the Weyl groups for \bide \, and \bidet;

\epw, \epv \, -- the determinants of the Weyl reflections $w$ and
$v$;

$C$, \CCt \, -- the Weyl chambers dominant with respect to $S$
and \SSt;

\ovC, \ovCCt \, -- the closures of the corresponding Weyl chambers;

$P_g$, $P_{\ggt}$ -- the weight lattices for $g$ and \ggt;

$\cal{E}$, $\widetilde{\cal{E}}$ -- the formal algebras associated to
 the weight lattices, $\cal{E}$ and $\widetilde{\cal{E}}$ are
 generated by the elements $e^{\beta}$, where $\beta$ is the
 fundamental weight, and the composition $e^{\beta} \cdot e^{\gamma}
 = e^{\beta + \gamma}$;

 \carla \,  -- the formal character of the representation \lla;

 \psil, \psitm \, -- the elements of formal algebras associated
 to the sets of anomalous weights for representations
 $\lla$, $\widetilde{L^{\mu}}$:
\begin{equation}
\psil = \sum_{w \in W} \epw
e^{w(\lambda + \rho) - \rho},   \label{e4}
\end{equation}
\begin{equation}
\psitm = \sum_{v \in V} \epv e^{v(\mu + \rrt) - \rrt)}.  \label{e5}
\end{equation}

For roots and weights of simple Lie algebras we use the standard
$e$-basis \cite{BUR}.

\section{Recursion Relations for Regular Injections}
\subsection{}
The initial decomposition (\ref{e1}) can be rewritten in terms of
formal characters \cite{NNN}
\begin{equation}
\carla = \sum_{\mu} n_{\mu} \carmu .
\label{e6}
\end{equation}
Applying the Weyl formula \cite{Kac}
\begin{equation}
\carks = \frac{\Psi^{\xi}}{ \proal }
\label{e7}
\end{equation}
and taking into account the injection
\bidet $\rightarrow$ \bide \, one gets the relation between the
anomalous elements \psil \, and \psitm \,
\begin{equation}
(\promi)^{-1} \psil = \sum_{\mu} n_{\mu} \psitm.
\label{e8}
\end{equation}
Using the basis $\{ e^{\xi} \}_{\xi \in P_{g}}$
of algebra $\cal{E}$ one can
expand both sides of the relation (\ref{e8}).
In the
left-hand side of the expansion of the first factor
\begin{equation}
(\promi )^{-1} = \sum_{\xi} K_{\ggt \subset g} (\xi) e^{-\xi}
\label{e10}
\end{equation}
gives rise to the so-called
Kostant-Heckman partition function \cite{Heck}.
The expression (\ref{e10})
toge\-ther with the formulas (\ref{e4}) and
(\ref{e5}) gives the desired
expansion of (\ref{e8}). Now consider only the
weights $\mu$ from
the dominant Weyl chamber \CCt, that is the anomalous weights
with $v=e$. Comparing the coefficients one gets the expression
for the relative multiplicity $n_{\mu}$ in terms of partition
function \Kgtg:
\begin{equation}
n_{\mu} = \sum_{W} \epw \Kgtg (w(\lambda + \rho) -
(\rho + \mu)).
\label{e11}
\end{equation}
Note that the relation (\ref{e8}) is valid on the whole weight
lattice $P_{g}$. Thus one can rewrite it in the form:
\begin{equation}
\sum_{\mu} n_{\mu} \psitm = \sum_{\mu} n_{\mu} \sum_{v \in V}
\epv e^{v(\mu + \rrt) - \rrt}
= \sum_{\xi} \widetilde{n_{\xi}} e^{\xi}.
\label{e12}
\end{equation}
Since all the weights $\{ v(\mu + \rrt) - \rrt \}  $ are
different the coefficients
$\nxti$ here are just the anomalous relative
multiplicities (see (\ref{e3})).
The relation (\ref{e12}) together with (\ref{e10})
and (\ref{e4}) shows that the expression (\ref{e11})
is true in all points of $P_g$ when $n$ is changed by
$\widetilde{n}$:
\begin{equation}
\nxti = \sum_{W} \epw \Kgtg (w(\lambda + \rho) -
(\rho + \xi)).
\label{e13}
\end{equation}
We would use this formula to
construct the recursion relation for \nxti.
First consider the expression
(\ref{e13}) for $\lambda = 0$. In this case $\nxti$
can be written explicitly as
the multiplicity of anomalous weight diagram for
the trivial subrepresentation $\widetilde{L^{0}}$:
\begin{equation}
\nxti = \sum_{V} \epv \delta_{\nu, \; v \rrt -\rrt}
= \sum_{W} \epw \Kgtg (w\rho  - (\rho + \xi)).
\label{e14}
\end{equation}
 One can extract the trivial term (with $w=e$)
and obtain the recursion relation
for the partition function \Kgtg
\begin{equation}
\Kgtg (\xi) = -\sum_{W \setminus e} \Kgtg (\xi + (w-1)\rho)
              + \sum_{V} \epv \delta_{\xi, \rho - v \rho}.
\label{e15}
\end{equation}
Returning to the expression (\ref{e13}) one easily comes to
the conclusion that the relation (\ref{e15}) induces the
recursion relation for anomalous relative multiplicities
\begin{equation}
\nxti = - \sum_{W \setminus e} \epw
\widetilde{n}_{\xi + (1-w)\rho} +
 \sum_{W,V} \epw \epv \delta_{\xi+(1-v)\rrt, \;
 w(\lambda + \rho) - \rho}.
\label{e16}
\end{equation}
This relation can be used in the explicit calculations of
multiplicities
$\nnti$ and $n_{\mu}$ and in some cases it is
effective. But the necessity
to estimate at each step the full set
of probe weights $\{ \xi + (1-w) \rho \}$
can make the whole process quite cumbersome.
Consider once more the formula (\ref{e8}).
The previous derivation was based on the properties
of the operator $(\promi)^{-1} $.
Now we shall use its inverse and taking into account
(\ref{e12}) rewrite the relation (\ref{e8}) in the following form
\begin{equation}
\psil = \sum_{w \in W} \epw
e^{w(\lambda + \rho) - \rho} = \promi \cdot \sum_{\xi}
\nxti e^{\xi}.
\label{e17}
\end{equation}
The first factor in (\ref{e17}) defines the finite set
 of weights $\Gamma ( g \subset \ggt ) $
whose structure depends only on the injection $ \ggt \subset g $,
\begin{equation}
\promi = 1 - \sum_{\gamma \in \Gamma } \mbox{sign} (\gamma )
e^{-\gamma}.
\label{e18}
\end{equation}
In these terms the formula (\ref{e17}) leads to the
following recursion relation
\begin{equation}
\nnti =
\sum_{\gamma \in \Gamma } \mbox{sign}(\gamma )
\widetilde{n}_{\nu + \gamma } + \sum_{w \in W} \epw
\delta_{\nu , \; w(\lambda + \rho ) - \rho }.
\label{e19}
\end{equation}
Its efficiency mainly depends on the possibility to construct
explicitly the set $\Gamma ( g \subset \ggt )$.

\subsection{}
To reveal the structure of $\Gamma$ let us use the denominator
identity \cite{Kac},
\begin{equation}
\proal = \Psi^{0}
\label{e20}
\end{equation}
to transform the expression (\ref{e18}):
\begin{equation}
\promi = 1 - \sum_{\gamma \in \Gamma } \mbox{sign}(\gamma )
e^{-\gamma} = (\prod_{\alpha \in \widetilde{\Delta}}
(1- e^{-\alpha}))^{-1} \cdot \Psi^{0}
\label{e21}
\end{equation}
The anomalous element $\Psi^{0}$ is $W$-invariant and can be
factorized with respect to $V \subset W$,
\begin{equation}
\Psi^{0} = \sum_{x \in X} \sum_{v \in V} \epsilon(v \cdot x)
e^{(v \cdot x - 1)\rho }.
\label{e22}
\end{equation}
Here $X$ is the factorspace $W/V$. This allows to present $\Gamma$
as the set of weight diagrams of representations $\widetilde{L}$:
\begin{equation}
1 - \sum_{\gamma \in \Gamma } \mbox{sign}(\gamma ) e^{-\gamma }
= \sum_{x \in X} \epsilon(x) (\prod_{\alpha \in \widetilde{\Delta}}
(1- e^{-\alpha}))^{-1} \sum_{v \in V} \epv e^{(vx\rho -
\rrt + (\rrt - \rho )) } = e^{\rrt - \rho} \sum_{x \in X}
\epsilon(x) ch \rule{1mm}{0mm} \widetilde{L}^{(x\rho - \rrt)}
\label{e23}
\end{equation}
The element $\Psi^{0}$ multiplied by $(\proal)^{-1}$ gives the weight
of the trivial representation $L^{0}$ of $g$, while the same
$\Psi^{0}$ multiplied by $(\prod_{\alpha \in \widetilde{\Delta}}
(1- e^{-\alpha}))^{-1}$ generates the assembly
$\Xi (\ggt \subset g)$ of weight diagrams
for representations of $\ggt$. To construct $\Xi (\ggt \subset g)$
one can use the auxiliary set $\Omega (\ggt \subset g)$
$$
\Omega = \{ 0, \alpha_{i_{1}}, \alpha_{i_{1}} + \alpha_{i_{2}},
\ldots , \alpha_{i_{1}} + \ldots + \alpha_{i_{m}} \mid
\alpha_{i_{l}} \in \Delta \setminus \bidet , m =
\mbox{card}(\Delta \setminus
\bidet ) \}.
$$
Fix the subset $\Omega'$ of dominant weights
$$
\Omega' = \{ \omega \in \Omega \mid \omega \in \ovCCt \}.
$$
Equip every $\omega \in \Omega $ with the sign
$\delta (\omega ) = \delta (\alpha_{i_{1}} + \ldots +
\alpha_{i_{k}}) = (-1)^{k+1} $.
Reduce $\Omega'$ to $\Omega_{r}'$ cancelling every
pair of weights in $ \Omega'$ that has opposite signs.
The subset of maximal weights
in $\Omega_{r}'$ is just the desired assembly of representations
$\Xi (\ggt \subset g)$ presented by their highest weights.
Obviously each $\xi \in \Xi (\ggt \subset g)$ has the form $\xi =
x\rho - \rho $ and it is easily seen that $\delta (\xi) =
\epsilon(x)$.

Let $\Phi^{\xi}$ be the weight diagram of the irrep
$\widetilde{L^{\xi}}$
with the highest weight $\xi \in \Xi (\ggt \subset g)$.
Due to the relation (\ref{e23}) the
set $\Gamma (\ggt \subset g)$ is obtained as the union of diagrams
$\Phi^{\xi}$
\begin{equation}
\Gamma = (\bigcup_{\xi \in \Xi (\ggt \subset g)} \Phi^{\xi})
\setminus \{ 0 \} .
\label{e24}
\end{equation}
Thus to find $\Gamma (\ggt \subset g)$ it is sufficient to construct
$\Xi (\ggt \subset g)$.
As we have seen the latter depends only on the structure of the
factorspace $W/V$ and the weights $\rho$ and $\rrt$.

Now we shall illustrate the situation treating different types of
regular maximal injections. For the injection $A_{n-1} \oplus u(1)
\rightarrow A_n$ the dimension of $W/V$ gives
$\mbox{card} (\Xi) = (n+1)!/n! = n+1$.
In the cases $A_{n-1} \oplus u(1) \rightarrow B_n , C_n , D_n$
the number of weights in $\Xi$ is proportional to
powers of "2":
$$
\mbox{card}(\Xi) = \left\{ \begin{array}{ll}
        2^{n} & \mbox{for }A_{n-1} \oplus u(1) \rightarrow B_n , \\
        2^{n} & \mbox{for }A_{n-1} \oplus u(1) \rightarrow C_n , \\
        2^{n-1} & \mbox{for }A_{n-1} \oplus u(1) \rightarrow D_n .
                           \end{array}
                   \right.
$$

The elementary analysis of the orbits of the Weyl group $V$
on the space of faithful representation of $W$ generated by
the weight $\rho$ leads to the following results.
\newtheorem{lemma}{Lemma}
\begin{lemma}
Put $\alpha_{0} \equiv (0, \ldots ,0) $ and let
$\{ \alpha_{1}, \ldots ,\alpha_{s} \} $ be the
ordered sequences of roots:
$$
\begin{array}{cl}
\{ e_{1}-e_{n+1}, e_{2}-e_{n+1}, \ldots , e_{n}-e_{n+1} \} &
\mbox{for} \;\; \Delta_{A_{n}} \setminus \Delta_{A_{n-1}} ,\\
\{ e_{1}+e_{2}, e_{1}+e_{3}, \ldots , e_{1}+e_{n}, e_{1}, & \\
e_{1}-e_{2}, e_{1}-e_{3}, \ldots , e_{1}-e_{n} \} &
\mbox{for} \;\; \Delta_{B_{n}} \setminus \Delta_{B_{n-1}} ,\\
\{ e_{1}+e_{2}, e_{1}+e_{3}, \ldots , e_{1}+e_{n}, 2e_{1}, & \\
e_{1}-e_{2}, e_{1}-e_{3}, \ldots , e_{1}-e_{n} \} &
\mbox{for} \;\; \Delta_{C_{n}} \setminus \Delta_{C_{n-1}} ,\\
\{ e_{1}+e_{2}, e_{1}+e_{3}, \ldots , e_{1}+e_{n}, & \\
e_{1}-e_{2}, e_{1}-e_{3}, \ldots , e_{1}-e_{n} \} &
\mbox{for} \;\; \Delta_{D_{n}} \setminus \Delta_{D_{n-1}} ;\\
\end{array}
$$
then the set $\Xi$ for $A_{n-1} \oplus u(1) \rightarrow A_n$
contains the weights
$$
\xi_{k} = \sum_{j=0}^{k} \alpha_{j} ; \; \; \; \; k=0, \ldots ,n
$$
while for $B_{n-1} \oplus u(1) \rightarrow B_n$ and
$C_{n-1} \oplus u(1) \rightarrow C_n$
$$
\xi_{k} = \sum_{j=0}^{k} \alpha_{j} ; \; \; \; \; k=0, \ldots ,2n-1
$$
and for $D_{n-1} \oplus u(1) \rightarrow D_n$
$$
\xi_{k} = \sum_{j=0}^{k} \alpha_{j} ; \; \; \; \; k=0, \ldots ,2n-2 $$
$$
\xi_{2n-1} = (n-1,1, \ldots ,1,-1).
$$
\end{lemma}
\begin{lemma}
In the sets $\Xi(A_{n-1} \oplus u(1) \rightarrow B_n, C_n, D_n) $
the zeroth weight is trivial
$$ \xi_{0} = (0, \ldots ,0);$$
the first weight has the form
$$ \xi_{1} = (p_{1}, p_{2}, \ldots ,p_{n}) =
\left\{ \begin{array}{ll}
(1,0,0,\ldots,0) & \mbox{for} \;\;
A_{n-1} \oplus u(1) \rightarrow B_n ,\\
(2,0,0,\ldots,0) & \mbox{for} \;\;
A_{n-1} \oplus u(1) \rightarrow C_n ,\\
(1,1,0,\ldots,0) & \mbox{for} \;\;
A_{n-1} \oplus u(1) \rightarrow D_n ;
       \end{array}
\right.   $$
while the others can be obtained from the following
relations:
$$ \xi_{2^{k}} = (\underbrace{p_{1}+k,1,1,\ldots,1}_{k+l},0,\ldots,0),
\;\;\; l = \left\{ \begin{array}{ll}
 1 & \left\{ \begin{array}{l}
           \mbox{for} \;\; A_{n-1} \oplus u(1) \rightarrow B_n  \\
           \mbox{for} \;\; A_{n-1} \oplus u(1) \rightarrow C_n
   \end{array}
   \right.    \\
 2 & \mbox{for} \;\; A_{n-1} \oplus u(1) \rightarrow D_n
             \end{array}
     \right. ,
$$
$$ \xi_{m}=\xi_{2^{k}+i}=\xi_{2^{k}}+(\xi_{i})_{shift} ; $$
$$ i=1,\ldots,
2^{k}-1 ;$$
where for every $\xi_{i}=(q_{1},q_{2},\ldots,q_{n})$ the
shifted weight $(\xi_{i})_{shift}$ is defined with the
coordinates
$$(\xi_{i})_{shift} = (0,q_{1},q_{2},\ldots,q_{n-1}) .$$
\end{lemma}

Note that in the framework of rules described in Lemma 1 and Lemma 2
the sets $\Xi$ are totally defined by the first weight $\xi_{1}$.

Now let us accumulate the information about $\Xi$'s for the
described types of injections in the following tables.

\newpage
\renewcommand{\arraystretch}{0.8}
\begin{tabular}{||c|c|c|c||} \hline
\raisebox{-2mm}{$\ggt \rightarrow g$}
& \multicolumn{2}{c|}{The
                  set $\Xi (\ggt \subset g)$ in terms of}
       & \raisebox{-2mm}{sign $\gamma(\xi)$}            \\
\cline{2-3}
       & highest weights $\xi$  &  Dynkin indices of
       $ \widetilde{L^{\xi}}$   &           \\
\hline
       & $ (0,0, \ldots ,0)$ & $([0,0, \ldots , 0],0)$ &        \\
$ A_{n-1} \oplus u(1)$
       & $ (1,0, \ldots ,0,-1)$ & $([1,0, \ldots , 0],1)$ & $ + $   \\
$ \rule{0.8cm}{0cm} \rightarrow A_{n} $
       & $ (1,1,0, \ldots ,0,-2)$ & $([0,1, \ldots , 0],2)$ & $ - $ \\
       & $\cdots$ & $\cdots$ & $\cdots$ \\
       & $ (1,1, \ldots ,1,-n)$ &
       $([0,0, \ldots , 0],n)$ & $ (-1)^{n+1} $ \\
\hline
       & $ (0,0, \ldots ,0)$ & $([0,0, \ldots , 0],0)$ &      \\
       & $ (1,0, \ldots ,0)$ & $([1,0, \ldots , 0],1)$ & $ + $      \\
       & $ (2,1,0, \ldots ,0)$ & $([1,1,0, \ldots , 0],3)$ & $ - $  \\
       & $ (2,2,0, \ldots ,0)$ & $([0,2,0, \ldots , 0],4)$ & $ + $  \\
$ A_{n-1} \oplus u(1)$
       & $ (3,1,1,0, \ldots ,0)$ &
       $([2,0,1,0, \ldots , 0],5)$ & $ + $      \\
$ \rule{0.8cm}{0cm} \rightarrow B_{n} $
       & $ (3,2,1,0, \ldots ,0)$ &
       $([1,1,1,0, \ldots , 0],6)$ & $ - $      \\
       & $ (3,3,2,0, \ldots ,0)$ &
       $([0,1,2,0, \ldots , 0],8)$ & $ + $      \\
       & $ (3,3,3,0, \ldots ,0)$ &
       $([0,0,3,0, \ldots , 0],9)$ & $ - $      \\
       & $\cdots$ & $\cdots$ & $\cdots$ \\
       & $ (n,n, \ldots ,n,n-1)$ & $([0,0, \ldots , 0,1],n^2-1)$
       & $ (-1)^{1/2(n^2+n)} $      \\
       & $ (n,n, \ldots ,n)$ & $([0,0, \ldots , 0],n^2)$
       & $ (-1)^{1/2(n^2+n+2)} $      \\
\hline
       & $ (0,0, \ldots ,0)$ & $([0,0, \ldots , 0],0)$ &      \\
       & $ (2,0, \ldots ,0)$ & $([2,0, \ldots , 0],2)$ & $ + $   \\
       & $ (3,1,0, \ldots ,0)$ & $([2,1,0, \ldots , 0],4)$ & $ - $  \\
       & $ (3,3,0, \ldots ,0)$ & $([0,3,0, \ldots , 0],6)$ & $ + $  \\
$ A_{n-1} \oplus u(1)$
       & $ (4,1,1,0, \ldots ,0)$ &
       $([3,0,1,0, \ldots , 0],6)$ & $ + $  \\
$ \rule{0.8cm}{0cm} \rightarrow C_{n} $
       & $ (4,3,1,0, \ldots ,0)$ &
       $([1,2,1,0, \ldots , 0],8)$ & $ - $  \\
       & $ (4,4,2,0, \ldots ,0)$ &
       $([0,2,2,0, \ldots , 0],10)$ & $ + $ \\
       & $ (4,4,4,0, \ldots ,0)$ &
       $([0,0,4,0, \ldots , 0],12)$ & $ - $ \\
       & $\cdots$ & $\cdots$ & $\cdots$ \\
       & $ (n+1, \ldots ,n+1,n-1)$ & $([0,0, \ldots , 2],n^2+n-2)$
       & $ (-1)^{1/2(n^2+n)} $      \\
       & $ (n+1, \ldots ,n+1)$ & $([0,0, \ldots , 0],n^2+n)$
       & $ (-1)^{1/2(n^2+n+2)} $      \\
\hline
       & $ (0,0, \ldots ,0)$ & $([0,0, \ldots , 0],0)$ &        \\
       & $ (1,1,0, \ldots ,0)$ &
       $([0,1,0, \ldots , 0],2)$ & $ + $      \\
       & $ (2,1,1,0, \ldots ,0)$ &
       $([1,0,1,0, \ldots , 0],4)$ & $ - $  \\
       & $ (2,2,2,0, \ldots ,0)$ &
       $([0,0,2,0, \ldots , 0],6)$ & $ + $ \\
$ A_{n-1} \oplus u(1)$
       & $ (3,1,1,1,0, \ldots ,0)$ &
       $([2,0,0,1,0, \ldots , 0],6)$ & $ + $ \\
$ \rule{0.8cm}{0cm} \rightarrow D_{n} $
       & $ (3,2,2,1,0, \ldots ,0)$ &
       $([1,0,1,1,0, \ldots , 0],8)$ & $ - $ \\
       & $ (3,3,2,2,0, \ldots ,0)$ &
       $([0,1,0,2,0, \ldots , 0],10)$ & $ + $ \\
       & $ (3,3,3,3,0, \ldots ,0)$ &
       $([0,0,0,3,0, \ldots , 0],12)$ & $ - $ \\
       & $\cdots$ & $\cdots$ & $\cdots$ \\
       & $ (n \! - \! 1, \! \ldots \! ,n
       \! -\! 1,n \! -\! 2,\! n \! -\! 2)$
       & $([0,\! \ldots \! ,0,1,0],n^2 \! - \! n \! - \! 2)$
       & $ (-1)^{1/2(n^2-n)} $ \\
       & $ (n-1, \ldots ,n-1)$ & $([0, \ldots , 0],n^2-n)$
       & $ (-1)^{1/2(n^2-n+2)} $ \\
\hline
\end{tabular}

\begin{center}
{\bf Table 1}
\end{center}
\newpage

\renewcommand{\arraystretch}{0.9}
\begin{tabular}{||c|c|c|c||} \hline
\raisebox{-2mm}{$\ggt \rightarrow g$}
& \multicolumn{2}{c|}{The
                  set $\Xi (\ggt \subset g)$ in terms of}
       & \raisebox{-2mm}{sign $\gamma(\xi)$}            \\
\cline{2-3}
       & highest weights $\xi$  &  Dynkin indices of
       $ \widetilde{L^{\xi}}$   &           \\
\hline
       & $ (0,0, \ldots ,0)$ & $([0,0, \ldots , 0],0)$ &      \\
       & $ (1,1,0, \ldots ,0)$ & $([1,0, \ldots , 0],1)$ & $ + $   \\
       & $ (2,1,1,0, \ldots ,0)$ &
       $([0,1,0, \ldots , 0],2)$ & $ - $  \\
       & $ (3,1,1,1,0, \ldots ,0)$ &
       $([0,0,1,0, \ldots , 0],3)$ & $ + $  \\
       & $\cdots$ & $\cdots$ & $\cdots$ \\
$ B_{n-1} \oplus u(1)$
       & $ (n-2,1,1, \ldots ,1,0)$ & $([0,0, \ldots ,1,0],n-2)$
       & $ (-1)^{n-1} $      \\
$ \rule{0.8cm}{0cm} \rightarrow B_{n} $
       & $ (n-1,1,1, \ldots ,1)$ &
       $([0, \ldots ,0,2],n-1)$ & $ (-1)^n $      \\
       & $ (n,1,1, \ldots ,1)$ &
       $([0, \ldots ,0,2],n)$ & $ (-1)^{n+1} $      \\
       & $ (n+1,1,1, \ldots ,1,0)$ & $([0, \ldots ,1,0],n+1)$
       & $ (-1)^{n+2} $      \\
       & $\cdots$ & $\cdots$ & $\cdots$ \\
       & $ (2n-2,1,0, \ldots ,0)$ & $([1,0, \ldots , 0],2n-2)$
       & $ - $      \\
       & $ (2n-1,0, \ldots ,0)$ & $([0,0, \ldots , 0],2n-1)$
       & $ + $      \\
\hline
       & $ (0,0, \ldots ,0)$ & $([0,0, \ldots , 0],0)$ &      \\
       & $ (1,1,0, \ldots ,0)$ & $([1,0, \ldots , 0],1)$ & $ + $   \\
       & $ (2,1,1,0, \ldots ,0)$ &
       $([0,1,0, \ldots , 0],2)$ & $ - $ \\
       & $ (3,1,1,1,0, \ldots ,0)$ &
       $([0,0,1,0, \ldots , 0],3)$ & $ + $      \\
$ C_{n-1} \oplus u(1)$
       & $\cdots$ & $\cdots$ & $\cdots$ \\
$ \rule{0.8cm}{0cm} \rightarrow C_{n} $
       & $ (n-1,1, \ldots ,1)$ &
       $([0, \ldots , 0,1],n-1)$ & $ (-1)^n $  \\
       & $ (n+1,1, \ldots ,1)$ &
       $([0, \ldots , 0,1],n+1)$ & $ (-1)^{n+1} $  \\
       & $ (n+2,1, \ldots ,1,0)$ &
       $([0, \ldots ,0,1,0],n+2)$
       & $ (-1)^{n+2} $  \\
       & $\cdots$ & $\cdots$ & $\cdots$ \\
       & $ (2n-1,1,0, \ldots ,0)$ &
       $([1,0, \ldots , 0],2n-1)$ & $ - $  \\
       & $ (2n,0, \ldots ,0)$ & $([0, \ldots , 0],2n)$ & $ + $ \\
\hline
       & $ (0,0, \ldots ,0)$ & $([0,0, \ldots , 0],0)$ &        \\
       & $ (1,1,0, \ldots ,0)$ & $([1,0, \ldots , 0],1)$ & $ + $  \\
       & $ (2,1,1,0, \ldots ,0)$ &
       $([0,1,0, \ldots , 0],2)$ & $ - $  \\
       & $\cdots$ & $\cdots$ & $\cdots$ \\
$ D_{n-1} \oplus u(1)$
       & $ (n-3,1, \ldots ,1,0,0)$ & $([0, \ldots ,0,1,0,0],n-3)$
       & $ (-1)^{n-2} $ \\
$ \rule{0.8cm}{0cm} \rightarrow D_{n} $
       & $ (n-2,1, \ldots ,1,0)$ & $([0, \ldots , 0,1,1],n-2)$
       & $ (-1)^{n-1} $ \\
       & $ (n-1,1, \ldots ,1,-1)$ & $([0, \ldots , 0,2,0],n-1)$
       & $ (-1)^n $ \\
       & $ (n-1,1, \ldots ,1,1)$ & $([0, \ldots , 0,2],n-1)$
       & $ (-1)^n $ \\
       & $ (n,1, \ldots ,1,0)$ & $([0, \ldots , 0,1,1],n)$
       & $ (-1)^{n+1} $ \\
       & $ (n+1,1, \ldots ,1,0,0)$ & $([0, \ldots ,1,0,0],n+1)$
       & $ (-1)^{n+2} $ \\
       & $\cdots$ & $\cdots$ & $\cdots$ \\
       & $ (2n-3,1,0, \ldots ,0)$ & $([1,0, \ldots ,0],2n-3)$
       & $ + $ \\
       & $ (2n-2,0, \ldots ,0)$ & $([0, \ldots ,0],2n-2)$
       & $ - $ \\
\hline
\end{tabular}

\begin{center}
{\bf Table 2}
\end{center}
\renewcommand{\arraystretch}{1.0}

\newpage

To conclude the general exposition of the recursion properties
of $\nmti$ we show the interdependence of the two recursion
formulas (\ref{e16}) and (\ref{e19}).
\begin{lemma}
The recursion relation (\ref{e16}) can be factorized with respect
to the subgroup $V$ of $W$ so that the summation over the
factorspace $W \setminus V$ is replaced by the summation over
$\Gamma$.
\end{lemma}
{\bf Proof.} Use the formula (\ref{e19}) to write down the
recursion relation for expression  \\
 $\sum_{v \in V} \epw
\widetilde{n}_{\nu + (1-v)\widetilde{\rho}} $
as a whole and extract the first term corresponding to
$v=e$:
\begin{eqnarray}
\nnti & = &
 - \sum_{v \in V, v \neq e} \epw
\widetilde{n}_{\nu + (1-v)\widetilde{\rho}}
- \sum_{v \in V} \sum_{\gamma \in \Gamma } \epv \;
\mbox{sign}(\gamma )
\widetilde{n}_{\nu + (1-v) \widetilde{\rho} +
\gamma }  \nonumber \\
 & & + \sum_{w \in W,v \in V} \epw \epv \delta_{\nu+(1-v)\rrt, \;
 w(\lambda + \rho) - \rho} \nonumber \\
 & = &
 \sum_{ V, \Gamma \cup \{ 0 \}} \epv \; \mbox{sign}(\gamma )
\widetilde{n}_{\nu + (1-v) \widetilde{\rho} + \gamma } +
 \nonumber \\
 & & + \sum_{w \in W,v \in V} \epw \epv \delta_{\nu+(1-v)\rrt, \;
 w(\lambda + \rho) - \rho} .
\label{el3}
\end{eqnarray}
 Comparing this expression with (\ref{e16}) we see that introducing
 $\Gamma$ one provides the factorization in the first term of the
 relation (\ref{e16}). Thus the relation (\ref{e19}) can be called
 {\em the factorized recursion formula} for anomalous relative
 multiplicities.

 \subsection{}

Both formulas (\ref{e16}) and (\ref{e19}) provide the
effective tools to treat the branching rules
decompositions for maximal regular injections. One can easily
estimate the relative capacities of these relations for
different pairs of
$g$ and $\ggt$. The result is that there are five families of
injections mostly favourable for the relation (\ref{e1}):
$A_{n-1} \oplus u(1)
\rightarrow A_n$, $B_{n-1} \oplus u(1)
\rightarrow B_n $, $C_{n-1} \oplus u(1)
\rightarrow C_n$, $D_{n-1} \oplus u(1)
\rightarrow D_n $, $A_{n-1} \oplus u(1)
\rightarrow B_n$. For these
five types the efficiency of the formula (\ref{e19}) increases with
growth of $n$ in comparison with that of (\ref{e16}). For the first
four types the ordinary decomposition methods are suitable (Gelfand-
Zeitlin procedure \cite{G-Z}, for example). So we shall concentrate
our attention on the last family: $A_{n-1} \oplus u(1)
\rightarrow B_n$.

To show the application of the factorized formula (\ref{e19}) in
details we shall start with the quite simple example. Consider the
injection  $ A_{1} \oplus u(1) \rightarrow B_{2} $. Fix the basic
roots of $B_{2}$
$$
S(B_{2}) = \{ \alpha_{1} = e_{1} - e_{2} , \alpha_{2} = e_{2} \}
$$
and the fundamental weights
$$
\{ \beta_{1} = e_{1} , \beta_{2} = \frac{1}{2} (e_{1} + e_{2}) \}.
$$
According to table 1 we have four highest weights in the set
 $\Xi ( A_{1} \oplus u(1) \rightarrow B_{2} )$:
$$
\Xi ( A_{1} \oplus u(1) \rightarrow B_{2} ) =
\{ (0,0),(1,0),(2,1),(2,2) \}
$$
Thus the set $\Gamma ( A_{1} \oplus u(1) \rightarrow B_{2} )$ contains
the weight diagrams of $A_{1} \oplus
u(1)$-representations $\Phi^{\xi_{1}} = ([1],1)$, $\Phi^{\xi_{2}} =
 ([1],3)$, $\Phi^{\xi_{s}} = ([0],4)$. ( The
$u(1)$-generator is normalized to have integer eigenvalues.)
$$
\Gamma = (\bigcup_{\xi \in \Xi} \Phi^{\xi}) \setminus \{ 0 \} =
\{ \gamma^{(sign \gamma)} \} =
\{ (1,0)^{(+)}, (0,1)^{(+)}, (2,1)^{(-)}, (1,2)^{(-)}, (2,2)^{(+)} \}.
$$
To simplify the following steps it is convenient to perform
further splitting of
the diagram of anomalous weights for subrepresentations
$\widetilde{L^{\mu}}$. This splitting is not unique,
one can choose an arbitrary vector
$\varepsilon \in C$ and the projections $a(\kappa)$ are obtained as
the scalar products
$$
\langle (\kappa - \lambda ), - \varepsilon \rangle = a(\kappa )
$$
for every $\kappa$ from to the weight lattice $P_g$. The weight
$\kappa$ is said to belong to the level $a(\kappa)$.
The ordering for the components in $\psil$
thus induced guarantees an inambiguous level by level
application of the recursion formula (\ref{e19}).
If $\Delta \setminus \bidet$ contains no positive
roots orthogonal to the boundary of $\ovC $, the auxiliary vector
$\varepsilon$ may be placed in the closure of $C$ as well.

Consider the irreducible representation $\lla$ of $B_{2}$ with
the highest weight $\lambda = (5/2, 1/2) $. The corresponding
anomalous weight diagram $\psil$ contains 8 vectors:
\begin{eqnarray*}
\lefteqn{\Psi^{(5/2, 1/2)} = \{ \psi^{\epw} \} = } \\
& & \{ (5/2, 1/2)^{(+)}, (-1/2, 7/2)^{(-)}, (5/2, -3/2)^{(-)},
(-5/2, 7/2)^{(+)}, (-1/2, -9/2)^{(+)}, \\ & & (-11/2, 1/2)^{(-)},
(-5/2, -9/2)^{(-)}, (-11/2, -3/2)^{(+)} \} .
\end{eqnarray*}
This is the case when splitting can be simplified. One can choose
$\varepsilon =(1,1) \in \ovC$ so
that $\langle \kappa ,\varepsilon \rangle$ becomes
proportional to the eigenvalues
of $u(1)$-generator in $\lla$
representation. Applying the factorized formula to obtain the
decomposition of $\lla$ we are interested in $\kappa$'s within the
Weyl chamber $\ovCCt$. Thus for our
example only the weights with nonnegative projection on
$\alpha_{1} = (1,-1)$ may have
positive multiplisities $n_{\kappa}$. On the zeroth level
the result is trivial
$$
\widetilde{n}_{(5/2,1/2)} = n_{(5/2,1/2)} = 1 ,$$
$$ \widetilde{n}_{(-1/2,7/2)} = -1 .
$$
On the next level (called the first) one finds two points in $\ovCCt$,
$(5/2,-1/2)$ and $(3/2,1/2)$,
where the formula (\ref{e19}) gives nonzero values for
$\widetilde{n_{\kappa}}$.
$$
n_{(5/2,-1/2)} = 1 ,\;\;\;\; n_{(3/2,1/2)} = 1.
$$
Similarly on the following two levels one finds
$$
n_{(1/2,1/2)} = 1 ,\;\;\;\; n_{(3/2,-1/2)} = 2
$$
and
$$
n_{(1/2,-1/2)} = 2 ,\;\;\;\; n_{(3/2,-3/2)} = 1.
$$
Due to the (reflection) symmetry of the weight diagram these four
levels give the sufficient information to write down the final
result:
\begin{eqnarray*}
\lefteqn{[2,1]_{\downarrow A_{1} \oplus u(1) } = } \\
 & & ([2],3) \oplus ([1],2) \oplus ([3],
2) \oplus  2([2],1) \oplus ([0],1) \oplus 2([1],0) \oplus ([3],0)
\\ & & \oplus 2([2],-1) \oplus ([0],-1) \oplus ([3],-2) \oplus
([1],-2) \oplus ([2],3).
\end{eqnarray*}
Here the numbers in the square brackets are the Dynkin indices and
the last term in the parenthesis is the eigenvalue of
$u(1)$-generator. Note that performing this reduction we need not
take into account the anomalous weights outside of the dominant
chamber $\ovCCt$. Such additional decomposition of recurrence
property takes place only when the vectors of the form
$ \{ \gamma + \xi \mid \xi \in \ovCCt, \gamma \in \Gamma \} $
do not reach the domain of anomalous weights
of $\ggt$ in $P_g \setminus (P_g \cap \ovCCt)$.

In \cite{L-F} an attempt was made to achieve the additional
decomposition of recurrence property in the situation when the
previous condition fails. The injection
$ A_3 \oplus u(1) \rightarrow D_4 $ was studied and
the recurrence relation connecting only relative multiplicities was
obtained. It is slightly different from that described by the formula
(\ref{e19}). But one faces great difficulties trying to obtain
such algorithm for other pairs of algebras.

In the appendix we bring the more complicated example considering on the
injection $ A_{3} \oplus u(1) \rightarrow B_{4} $. This demonstrates
the efficiency of the decomposition algorithm based on the factorized
recurrence formula (\ref{e19}). The whole computation is relatively
simple and can be easily computerized.
The nonmaximal regular injections and special injections can be
treated similarly. The detailed study of recurrence relations in
these cases will be presented in the forthcoming publication.

\underline{\bf Acknowledgments}

   Supported in part by Russian Foundation for
 Fundamental Research, Grant N 95-01-00569a . Supported in part by the
 International Science Foundation, Grant N U9J000. \\ \\

\newpage
\noindent {\large \bf Appendix. \\
The Injection $ A_{3} \oplus u(1)  \rightarrow B_{4} $.}

The positive root systems for this maximal regular injection in the
standard $e$-basis can be written as follows:
\begin{eqnarray*}
\Delta (B_{4}) & = & \{ e_{i}, e_{j}-e_{k}, e_{j}+e_{k} \} \\
\widetilde{\Delta}(A_{3}) & = & \{ e_{j}-e_{k} \} \\
\Delta \setminus \widetilde{\Delta} & = &
\{ e_{i}, e_{j}+e_{k} \}  \\
& & i, j, k = 1, \ldots ,4; \;\; j < k
\end{eqnarray*}
According to Lemma 2 the set $\Xi$ contains 16 irreducible
representations of $ A_{3} \oplus u(1) $ enumerated by their highest
weights $\xi_{i}$ $ (i = 0, \ldots ,15)$ :
\begin{eqnarray*}
\lefteqn{\Xi (A_{3} \oplus u(1) \rightarrow B_{4}) = }  \\
 & & \{ (0,0,0,0); (1,0,0,0); (2,1,0,0);  (2,2,0,0);
(3,1,1,0);  \\  & & (3,2,1,0);
(3,3,2,0); (3,3,3,0); (4,1,1,1);
(4,2,1,1);  \\  & & (4,3,2,1); (4,3,3,1);
(4,4,2,2); (4,4,3,2); (4,4,4,3);  \\  & & (4,4,4,4) \}.
\end{eqnarray*}
In the "fan" $\Gamma (A_{3} \oplus u(1) \rightarrow B_{4})$ the
weights $\gamma$ for each $\Phi^{\xi_{s}}$ bear the same sign,
thus the signs can be attributed to the
representations $\widetilde{L^{\xi_{s}}}$:
\begin{eqnarray*}
\lefteqn{\{ \widetilde{L^{\xi_{s}}};
 \mbox{sign}(\gamma) \}_{s=1,\ldots ,15} = }  \\
 & & \{ (([1,0,0],1);(+)), (([1,1,0],3);(-)), (([0,2,0],4);(+)),
(([2,0,1],5);(+)),  \\
 & & (([1,1,1],6);(-)), (([3,0,0],7);(-)), (([0,1,2],8);(+)),
(([2,1,0],8);(+)),  \\
 & & (([0,0,3],9);(-)), (([1,1,1],10);(-)), (([1,0,2],11);(+)),
(([0,2,0],12);(+)),  \\
 & &  (([0,1,1],13);(-)), (([0,0,1],15);(+)), (([0,0,0],16);(-))
\}.
\end{eqnarray*}
The nontrivial multiplicities of $\gamma \in \Gamma$ must also be
taken into account.
To obtain the splitting one can chose the vector
$\varepsilon = (1,1,1,1)$.
We write down
explicitly only those weights $\gamma$ that describe the first and
the second levels of decomposition:
\begin{eqnarray*}
\lefteqn{\Gamma (A_{3} \oplus u(1) \rightarrow B_{4}) =
\bigcup_{\xi_{s}, s = 1, \ldots ,15}
\Phi^{\xi_{s}; \; \mbox{sign}(\gamma)} = }
 \\ & & \{ (1,0,0,0)^{(+)}, (0,1,0,0)^{(+)}, (0,0,1,0)^{(+)},
(0,0,0,1)^{(+)}, (2,1,0,0)^{(-)},
 \\  & & (2,0,1,0)^{(-)},
(2,0,0,1)^{(-)}, (1,0,0,2)^{(-)}, (0,1,0,2)^{(-)}, (0,0,1,2)^{(-)},
 \\ & & (1,2,0,0)^{(-)}, (0,2,1,0)^{(-)},
(0,2,0,1)^{(-)}, (1,0,2,0)^{(-)}, (0,1,2,0)^{(-)},  \\  & &
(0,0,2,1)^{(-)}, 2(1,1,1,0)^{(-)}, 2(1,1,0,1)^{(-)},
2(1,0,1,1)^{(-)}, 2(0,1,1,1)^{(-)}, \ldots   \\  & &
\ldots ,  (4,4,4,4)^{(-)} \}.
\end{eqnarray*}
Consider for example the irrep $\lla(B_{4})$ with $\lambda = (2,2,1,
0)$; dim $\lla =(\underline{1650})$. It has 11 levels. Due to
the reflection symmetry ($(k_{1},k_{2},k_{3},k_{4})
\leftrightarrow (-k_{1},-k_{2},-k_{3},-k_{4})$) of the weight
diagram it is sufficient to study only 6 of them.

Contrary to the previous case the recurrence procedure cannot
be performed separately for the relative multiplicities
$n_{\kappa}$, that is for $\widetilde{n_{\kappa}}$ with
$\kappa \in \ovCCt $.
Nevertheless the
calculations involving the anomalous weights (in other Weyl chambers)
can be simplified considerably thanks to the following two
considerations:
\begin{itemize}
  \item In the level by level recursive procedure after the evaluation
  of $\widetilde{n}$ in $\ovCCt$ the anomalous relative
   weights in other points can be obtained using the
   Weyl group $V$.
  \item Using the exterior contour of $\Gamma$ one can easily fix the
   domain of the weight space that can contribute to the relative
   multiplicities $n_{\kappa}$ and
   pay no attention to the weights out of this domain.
\end{itemize}
For the injection $ A_{3} \oplus u(1) \rightarrow B_{4} $ only the
anomalous weights with nonnegative first coordinate must be taken into
account to obtain $\widetilde{n} \in
\ovCCt$. Thus for every obtained highest weight $\mu$ one must
find only twelve anomalous points of $\psim$ to be able
to proceed the recursion
on the next level. On the zeroth level besides the
highest weight
$$
n_{(2,2,1,0)} = \widetilde{n}_{(2,2,1,0)} = 1
$$
one must also calculate (using the Weyl group $V$)
the anomalous relative multiplicities:
$$
\widetilde{n}_{(2,2,-1,2)} = -1 ,\;\;
\widetilde{n}_{(2,0,3,0)} = -1 ,\;\;
\widetilde{n}_{(2,0,-1,4)} = 1 ,\;\;
\widetilde{n}_{(2,-2,3,2)} = 1 , $$
$$ \widetilde{n}_{(2,-2,1,4)} = -1 ,\;\;
\widetilde{n}_{(1,3,1,0)} = -1 ,\;\;
\widetilde{n}_{(1,3,-1,2)} = 1 ,\;\;
\widetilde{n}_{(1,0,4,0)} = 1 , $$
$$ \widetilde{n}_{(1,0,-1,5)} = -1 ,\;\;
\widetilde{n}_{(1,-2,4,2)} = -1 ,\;\;
\widetilde{n}_{(1,-2,1,5)} = 1 .
$$
After this the formula (\ref{e19}) can be directly applied to
obtain the anomalous relative multiplicities of the first level
and among them --
$$
n_{(2,2,0,0)} = 1 ,\;\; n_{(2,1,1,0)} = 1 ,
$$
The 24 anomalous points of these two representations
fix the decomposition on the second level:
$$
n_{(2,1,0,0)} = 2 ,\;\; n_{(1,1,1,0)} = 1 ,\;\; n_{(2,1,1,-1)} = 1 ,
\;\; n_{(2,2,0,-1)} = 1 ,
$$
and so on.
For example, on the third level for the weight $(1,1,0,0)$  the
formula (\ref{e19}) leads to the following relation:
$$
\widetilde{n}_{(1,1,0,0)} = n_{(1,1,0,0)} =
n_{(1,1,1,0)} + n_{(2,1,0,0)} - 2n_{(2,2,1,0)} -
\widetilde{n}_{(1,3,1,0)} = 2 .
$$
The final result is
\begin{eqnarray*}
\lefteqn{[0,1,1,0]_{\downarrow A_{3} \oplus u(1)} = } \\ & &
([0,1,1],5) \oplus ([0,2,0],4) \oplus ([1,0,1],4) \oplus
2([1,1,0],3)  \\ & & \oplus ([0,0,1],3) \oplus ([1,0,2],3) \oplus
([0,2,1],3) \oplus ([2,0,0],2)  \\ & & \oplus  2([0,1,0],2) \oplus
([0,0,2],2) \oplus 2([1,1,1],2) \oplus 2([1,0,0],1)  \\ & & \oplus
2([2,0,1],1) \oplus  3([0,1,1],1) \oplus ([1,1,2],1) \oplus
([1,2,0],1)  \\ & & \oplus ([0,0,0],0) \oplus 3([1,0,1],0) \oplus
2([0,2,0],0) \oplus ([2,0,2],0)  \\ & & \oplus ([2,1,0],0) \oplus
([0,1,2],0) \oplus 2([1,0,2],-1) \oplus  2([0,0,1],-1)  \\ & & \oplus
2([1,1,1],-2) \oplus ([2,0,0],-2) \oplus 2([0,1,0],-2) \oplus
([0,0,2],-2)  \\ & & \oplus ([1,2,0],-3) \oplus ([2,0,1],-3) \oplus
([1,0,0],-3) \oplus 2([0,1,1],-3)  \\ & & \oplus ([1,0,1],-4) \oplus
([0,2,0],-4) \oplus ([1,1,0],-5).
\end{eqnarray*}

\end{document}